\input amstex
\loadbold
\documentstyle{amsppt}
\TagsOnRight

\topmatter

\title
On finite-sheeted covering mappings onto solenoids
\endtitle

\rightheadtext{finite-sheeted coverings of solenoids}

\author
R. N. Gumerov
\endauthor

\address
Department \ of \ Mechanics \ and \ Mathematics, \ Kazan State University,
\linebreak Kremlevskaya 18, \ Kazan, \ 420008, \ Tatarstan, \ Russian Federation
\endaddress

\email
renat.gumerov\@ksu.ru
\endemail



\vskip 0.5 cm

\abstract
We study limit mappings from a solenoid onto itself.
It is shown that each equivalence class of finite-sheeted covering mappings
from connected topological spaces onto a solenoid is determined by a limit
mapping. Properties of periodic points of limit mappings are also studied.
\endabstract

\endtopmatter

\document

\head
Introduction
\endhead

There \, are \, various \, ways \, one \, can \, view \, solenoids  (\, see, e.g., \cite{4},
\cite{11, (10.12)}\,). A solenoid may be defined as follows.
Let $P=(\, p_1, p_2, \ldots)$ be a sequence of prime numbers (1 not being
included as a prime). The {\it $P$-adic solenoid} $\Sigma_P$ is the inverse
limit of the inverse sequence
$$
\CD
\Bbb S^1 @< f_1^2 << \Bbb S^1 @< f_2^3 << \Bbb S^1 @< f_3^4 << \cdots ,
\endCD
\tag {$1$}
$$
where $\Bbb S^1$ is the unit circle considered as a subspace of the space
$\Bbb C$ of all complex numbers endowed with the natural topology, and every
bonding mapping $f_n^{n+1}$ is given by $f_n^{n+1}(z)=z^{p_n}$ for each
$z \in \Bbb S^1$, $ n \in \Bbb N$. If $p_n=2$ for all $n\in \Bbb N$, then
the solenoid $\Sigma_P$ is said to be {\it dyadic}. As is well known,
the $P$-adic solenoid is a metric continuum which is not locally connected
at any point. Sequences $P$ and $Q$ of prime numbers are said to be
{\it equivalent} \ (written $P \sim Q$) if a finite number of terms can be
deleted from each sequence so that every prime number occurs the same number
of times in the deleted sequences. The solenoids $\Sigma_P$ and $\Sigma_Q$
are homeomorphic iff $P \sim Q$ \ \ \cite{2}, \cite{14} (see also \cite{1}).
The solenoid is a compact abelian group under the coordinatewise
multiplication with the identity $(1,1,\ldots)$. The condition $P \sim Q$
is also a criterion for existing of topological isomorphism between
the topological groups $\Sigma_P$ and $\Sigma_Q$ \cite{2}.

For given $k\in \Bbb N$ let us consider the limit mapping
$h_P^k:\Sigma_P \to \Sigma_P$ induced by the mapping $\{h_n^k: n\in \Bbb N\}$
between two copies of (1) :
$$
\CD
\Bbb S^1 @<f_1^2<< \Bbb S^1 @<f_2^3<< \Bbb S^1 @<f_3^4<< \cdots \ \ \Sigma_P \\
@Vh^k_1VV                  @VVh^k_2V            @VVh^k_3V                   \ \ \ \ \ \  @VVh_P^kV  \\
\Bbb S^1 @<f_1^2<< \Bbb S^1 @<f_2^3<< \Bbb S^1 @<f_3^4<< \cdots \ \ \Sigma_P ,
\endCD \tag {$2$}
$$
where \ $h^k_n$ \ is the $k$-th potency mapping for each $n\in \Bbb
N$, that is, \ $h^k_n(z)=z^k$ \ for every $z \in \Bbb S^1$.

This note deals with the limit mappings $h^k_P$. Note that, for each
$k \in \Bbb N$, the limit mapping \ $h^k_{(2,2,\ldots)}$ \ of the dyadic
solenoid is a covering mapping \cite{17}, \cite{3}.
Theorem 19 of \cite{3} states that the dyadic solenoid admits any odd
degree $k$ covering mapping of the form \ $h^k_{(2,2,\ldots)}$, while it
does not admit any covering mapping \ $h^k_{(2,2,\ldots)}$ \ of an even
degree $k$. For each $k\in \Bbb N$ the set of periodic points
of \ $h^k_{(2,2,\ldots)}$ \ is dense in \ $\Sigma_{(2, 2,\ldots)}$ \
\cite{17}, \cite{3}.

The purpose of this note is to extend the above results concerning
$h^k_{(2,2,\ldots)}$ to the limit mappings $h^k_P$ of an arbitrary
$P$-adic solenoid.

\head
1. Preliminaries
\endhead

In this section we establish some notation that is used throughout.
As usual, we denote by $\Bbb N$ the set of all positive integers.
For a complex number $z \in \Bbb C$ and $n\in \Bbb N$ we denote by \
$\root n \of{z}$ \ the set of all values of the $n$-th root of $z$.

It is well known that in studying of finite-sheeted covering mappings from
connected topological spaces onto solenoids there is no loss of generality in
assuming that covering spaces are metric continua, i.e., metric compact
connected spaces. So all topological spaces are assumed to be metric.
A {\it mapping} between two spaces always means a continuous function.
For basic notions of the theory of inverse limit spaces we refer the reader
to \cite{5, Chapter VIII} and \cite{6, Chapter 2}. Let $\{X_n, \pi_n^{n+1}\}$
\ and \ $\{Y_n, \rho_n^{n+1}\}$ \ be inverse sequences of spaces $X_n$ and $Y_n$
with bonding mappings
$\pi_n^{n+1}:X_{n+1}\to X_n$ and $\rho_n^{n+1}:Y_{n+1}\to Y_n$,
where $n\in \Bbb N$. We denote by $X_\infty$ and $Y_\infty$ the inverse
limit spaces of $\{X_n, \pi_n^{n+1}\}$ \ and \ $\{Y_n,\rho_n^{n+1}\}$ \
respectively.
Recall that a sequence $\{\sigma_n:X_n \to Y_n \mid n\in \Bbb N\}$ of mappings
is called \ a {\it mapping} \ between \ the inverse sequences
$\{X_n, \pi_n^{n+1}\}$ \ and \ $\{Y_n, \rho_n^{n+1}\}$ \ if
\ $\sigma_n\circ \pi_n^{n+1}=\rho_n^{n+1}\circ \sigma_{n+1}$ \ for all \ $n\in \Bbb N$. \
Then \ there \ exists \ a \ mapping \
$\sigma_\infty:X_\infty \to Y_\infty:(x_1,x_2,\ldots) \mapsto (\sigma_1(x_1),\sigma_2(x_2),\ldots)$.
It is called the {\it limit mapping} induced by $\{\sigma_n: n\in \Bbb N\}$.

For any $m,n\in \Bbb N$ satisfying $m\leq n$, we denote by $f_m^n$
the bonding mapping of the inverse sequence (1).
Thus, for each $n\in \Bbb N$, the mapping $f_n^n$ is the identity
on $\Bbb S^1$ and $f_l^n=f_l^m\circ f_m^n$ for all $l,m,n\in \Bbb N$
satisfying $l \leq m \leq n$.

The identity $(1,1,\ldots)$ of a solenoid is denoted by $e$.

Recall that a surjective mapping $f:X\to Y$ between spaces $X$ and $Y$ is
called:

--- $k$-{\it to}-$1$, where $k\in \Bbb N$, provided that $card \, f^{-1}(y)=k$
for each $y\in Y$;

--- a {\it finite-sheeted covering mapping} \ if \ it is a $k$-{\it sheeted} \
(\, $k$-{\it fold}\,) {\it covering mapping} for some $k\in \Bbb N$; that is,
every point $y\in Y$ has an open neighborhood $W$ in $Y$ such that the inverse
image $f^{-1}(W)$ can be written as the union of $k$ disjoint open subsets of
$X$ each of which is mapped homeomorphically onto $W$ under $f$.

The number $k$ in the above definitions is called a {\it degree} of
the mapping. Finite-sheeted covering mappings \ $f_1:X_1\to Y$ and
$f_2:X_2\to Y$ are said to be {\it isomorphic} (or {\it equivalent}) \ if
there is a homeomorphism \ $g:X_1 \to X_2$ \ such that \ $f_1= f_2 \circ g$.

Throughout this note $P=(\, p_1,p_2,\ldots)$ denotes a sequence of prime
numbers. We say that a prime number $p$ {\it occurs infinitely often} in
$P$ if $p=p_n$ for infinitely many terms $p_n$, $n \in \Bbb N$. We denote
by $S(P)$ a subset of $\Bbb N \smallsetminus \{1\}$ consisting of all prime
numbers which do not occur infinitely often in $P$. In other words,
$q\in S(P)$ iff $q > 1$ is a prime number such
that, for some $m\in \Bbb N$, we have $q\neq p_n$ for all $n\geq m$.

\head
2. Covering mappings.
\endhead

The following proposition is an analog of \cite{17, Lemma 1} and
\cite{3, Proposition 8}.

\proclaim{Proposition 1} For each $k\in \Bbb N$
the limit mapping \ $h^k_P: \Sigma_P \to \Sigma_P$ \
is a finite-sheeted covering mapping, and its degree is at most $k$.
\endproclaim
\demo{Proof} Let us fix $k\in \Bbb N$. Clearly, it suffices to
show that $h^k_P$ is an $m$-to-1 surjective open mapping for
some $m\leq k$ \ (cf. \cite{16, Chapter X, \S 6}, \cite{3, Proposition 1}).

Since $\Bbb S^1$ is compact and each mapping $h_n^k$ in (2) is surjective,
the limit mapping $h^k_P$ is surjective too \cite{6, Theorem 3.2.14}.
Because $h_P^k$ is a continuous homomorphism from a compact group onto itself,
$h_P^k$ is an open mapping \cite{11, (5.29)} and the equality
$card\, (h^k_P)^{-1}(y)=card\, (h^k_P)^{-1}(e)$ holds for each $y\in \Sigma_P$.
Note, if $(x_1,x_2,\ldots)\in (h^k_P)^{-1}(e)$,
then $x_n \in \root k \of{1}$ for all $n\in \Bbb N$. Using these
observations, one can easily see that $h_P^k$ is an $m$-to-1 mapping
and $m\leq k$.  \qed

\enddemo

We shall determine degrees of the limit mappings $h^k_P$. In order to do this
we find out cardinalities of fibers $(h^k_P)^{-1}(e)$, $k\in \Bbb N$.
(\, See also \cite{3, Statements 15, 17, 18}\,).

\proclaim{Proposition 2} If $k$ is a prime number which occurs infinitely
often in the sequence $P$, then the limit mapping \
$h^k_P: \Sigma_P \to \Sigma_P$ \ is a homeomorphism.
\endproclaim
\demo{Proof} We claim that the set $(h^k_P)^{-1}(e)$ consists of only one
point $e$. To show this we suppose that
$z=(z_1, z_2, \ldots) \in (h^k_P)^{-1}(e)$. Then $z_n^k=1$
for all $n\in \Bbb N$. By assumption, for given $n\in \Bbb N$,
there is an integer $m \geq n$ such that $f_m^{m+1}$ in (1) is
the $k$-th potency mapping. Therefore, we have
$$z_n=f_n^m(f_m^{m+1}(z_{m+1}))=f_n^m(z_{m+1}^k)=f_n^m(1)=1.$$
Since these equalities hold for each  $n\in \Bbb N$, it follows that $z=e$,
as claimed.

Thus $card\, (h^k_P)^{-1}(e)=1$ and, by Proposition 1, the limit mapping
$h^k_P$ is a homeomorphism. \qed
\enddemo

The verification of the following lemma is straightforward (cf. \cite{3,
Fact 16}).

\proclaim{Lemma} If \ $k=l\cdot m$ \ for some $l, m \in \Bbb N$, then
$h^k_P=h^l_P \circ h^m_P.$
\endproclaim

Combining Lemma and Proposition 2, we have :

\proclaim{Proposition 3} If each prime divisor of \ $k\in \Bbb N$ \
occurs infinitely often in the sequence $P$, then the limit
mapping \ $h^k_P: \Sigma_P \to \Sigma_P$ \ is a homeomorphism.
\endproclaim

Let $n\in \Bbb N$ and $\xi_n=\cos\frac{2\pi}{n}+i\sin\frac{2\pi}{n}$,
where $i$ is the imaginary unit, i.e., $i^2=-1$. The set $\root n \of{1}$ is
a multiplicative cyclic group $\{1, \xi_n, \xi_n^2, \ldots,\xi_n^{n-1}\}$
generated by $\xi_n$. For $m\in \Bbb N$, we define a homomorphism
$\psi_m:\root n \of{1} \ \to \ \root n \of{1}$ by setting
$\psi_m(\xi_n^j)=\xi_n^{jm}$, where $j \in \{1,2,\ldots,n\}.$
Let $n$ and $m$ be relatively prime. In this case, one can easily see
that $\psi_m$ is injective. This implies that $\psi_m$ is an
automorphism, i.e., a bijective homomorphism from a group onto
itself. Therefore, there exists an automorphism $\phi_m:\root n
\of{1} \ \to \ \root n \of{1}$, which is the inverse of $\psi_m$.
In other words, for each $a\in \root n \of{1}$ we have exactly
one element $b$ of the group $\root n \of{1}$ whose $m$-th power
is equal to $a$. The value of the automorphism $\phi_m$ at $a$ is
just the element $b$.

\proclaim{Proposition 4} If $k \in S(P)$, then the limit
mapping \ \ $h^k_P:\Sigma_P \to \Sigma_P$ \ is a $k$-fold covering mapping.
\endproclaim

\demo{Proof} Denote by $m\in \Bbb N$ the least number $l$ such
that $k\neq p_n$ for all $n\geq l$.

If $m=1$, i.e., $k$ is not a term of the sequence $P$,
then for each $p_n \in P$, the numbers $p_n$ and $k$ are relatively
prime. Let $\phi_{p_n}$ be the inverse of the automorphism
$\psi_{p_n}:\root k \of{1} \to \root k \of{1}: z \mapsto z^{p_n}$. We have
$(\phi_{p_n}(z))^{p_n}=z$ for every $z\in \Bbb S^1$.
It is easy to see that for each $k$-th root of unity
$\xi_k^j$, where $j \in \{1,2,\ldots,k\}$, the element
$$
(\xi_k^j,  \phi_{p_1}(\xi_k^j),  \phi_{p_2} \circ \phi_{p_1}(\xi_k^j), \ldots)
$$
lies in the set $(h_P^k)^{-1}(e)$. Clearly, in this way
one obtains $k$ distinct points of the fiber $(h_P^k)^{-1}(e)$.
Hence \ $card\, (h_P^k)^{-1}(e) \geq k$. On the other hand, by Proposition 1,
$card \, (h_P^k)^{-1}(e)\leq  k$ and the desired conclusion follows .

If $m > 1$, then we consider the sequence $Q=(p_m, p_{m+1}, \ldots)$
which is equivalent to $P$. As mentioned in Introduction, there exists
a topological isomorphism $\rho:\Sigma_P \to \Sigma_Q$ of
topological groups ( take $\rho$ to be the shift of the form
$\rho((x_1, x_2, \ldots))=
(x_{m},x_{m+1},\ldots)$, \ where \ $(x_1, x_2, \ldots)\in \Sigma_P $).
The inverse of $\rho$ is denoted by $\tau$.

One can readily verify that the following diagram

$$
\CD
\Sigma_P   @>\rho>>   \Sigma_Q \\
@Vh_P^kVV   @VVh_Q^kV \\
\Sigma_P @>\rho>> \Sigma_Q
\endCD \tag{$3$}
$$
is commutative, i.e. $\rho \circ h_P^k=h_Q^k \circ \rho$.
According to the first part of this proof, the limit mapping \ $h^k_Q$ \
is a $k$-fold covering mapping.

It follows from the commutativity of the diagram (3) that
$$
(h_P^k)^{-1}(e)=(h_P^k)^{-1}(\tau(e))=\tau((h_Q^k)^{-1}(e)).
$$
But the mapping \ $\tau$ \ is injective and $card \, (h_Q^k)^{-1}(e)=k$.
Thus the equality  \ \ \ $card \, (h_P^k)^{-1}(e)=k$ holds.
Therefore, by Proposition 1, $h^k_P$ is a $k$-fold covering mapping. \qed
\enddemo

As an immediate consequence of the above results we have the following
theorem.

\proclaim{Theorem 1} Let $k\in \Bbb N$ be given. The $P$-adic solenoid admits
a $k$-fold covering mapping of the form $h^k_P$ if and only if $k$ has
no a prime divisor which occurs infinitely often in the sequence $P$.
\endproclaim

It is interesting to restate Theorem 1. To do this we recall
some facts. We refer to \cite{11} for a basic material on topological groups.
The character group of $\Sigma_P$, consisting of all continuous homomorphisms
from $\Sigma_P$ into $\Bbb S^1$, is topologically isomorphic to the discrete
group $F_P$ of $P$-adic rationals (\, all rationals of the form
$m/(p_1p_2\ldots p_n)$, where $m$ is an integer and $n\in \Bbb N$ \, )
\cite{11, (25.3)}. It follows from the Pontrjagin duality that two
solenoids $\Sigma_P$ and $\Sigma_Q$ are topologically isomorphic if and only
if the additive groups $F_P$ and $F_Q$ are isomorphic.

An additive abelian group $G$ is said to be $n$-{\it divisible} provided
that for each element $g\in G$ there is an element $g' \in G$ such that
$ng^\prime=g$. By number-theoretic considerations one can see
that, for a prime number $p$, the group $F_P$ is
$p$-divisible if and only if $p$ occurs infinitely often in $P$.

\proclaim{Theorem $\boldkey 1^{\boldsymbol\prime} $ } Let $k\in \Bbb N$
be given. The $P$-adic solenoid admits a $k$-fold covering mapping
of the form $h^k_P$ if and only if $k$ has no a prime divisor $q$
such that the group of $P$-adic rationals is $q$-divisible.
\endproclaim

\remark{Remark 1} Let $k\geq 2$ and let $G$ be a compact connected abelian
group. Suppose that the character group of $G$ is $k$-divisible.
Using Theorem 1 of \cite{9}, it can be shown that there is no a $k$-sheeted covering
mapping from a connected Hausdorff topological space onto $G$ \cite{10}.
\endremark

In what follows {a \it finite-sheeted connected covering
of a solenoid} means a finite-sheeted covering mapping from a connected space
onto a solenoid. We conclude this section with a theorem which summarizes
the results concerning finite-sheeted connected coverings of solenoids.
For each $k\in \Bbb N$ that is not divisible by any of
the primes that occur infinitely often in the sequence $P$ there is just one
(up to equivalence) $k$-fold connected covering of $\Sigma_P$,
and these are the only finite-sheeted connected coverings of $\Sigma_P$
\ (\, see \cite{7, Example 2}, \cite{8, p. 82, Theorem 3} and
\cite{15, Proposition 2.2}\,).
This fact is a corollary of the theory of overlays
(\, We refer to \cite{7}, \cite{8} and \cite{12} for this theory\,).
According to the above results we have the following theorem
(with the same dichotomy of positive integers $k$ as in Theorems 1 and
$1^\prime $).

\proclaim{Theorem 2} Let $P$ be a sequence of prime numbers and
$k\in \Bbb N$. If $k$ is a multiple of some prime number which
occurs infinitely often in $P$, then there is no a $k$-fold
connected covering of the $P$-adic solenoid $\Sigma_P$.
Otherwise, the limit mapping $h^k_P$ is a $k$-fold connected covering and,
moreover, each $k$-fold connected covering of $\Sigma_P$ is equivalent
to $h^k_P$.
\endproclaim

\remark{Remark 2} The fact that each finite-sheeted connected
covering of $\Sigma_P$ is isomorphic to some limit mapping $h^k_P$
(and, as a consequence, Theorem 2) can be proved without using
the theory of overlays. In order to show this we make use of the approximate
construction sketched in \cite{9} (\, A detailed account of this construction
for a finite-sheeted covering mapping from a connected Hausdorff topological
space onto an arbitrary compact connected group is contained in \cite{10}\, ).
For the sake of completeness, we outline the proof of the above-mentioned fact
as follows.

Let $f:X\to \Sigma_P$ be a $k$-fold connected covering.
There exists an inverse sequence $\{X_n,g_n^{n+1}\}$ and a mapping
$\{g_n:X_n \to S^1 \mid n\in \Bbb N\}$ between $\{X_n,g_n^{n+1}\}$
and the inverse sequence (1) such that the properties listed below are fulfilled:
1) for each $n\in \Bbb N$ the space $X_n$ is connected and locally pathwise
connected and $g_n:X_n\to \Bbb S^1$ is a $k$-fold covering mapping;
2) the $k$-fold covering mapping $f$ is isomorphic to the limit mapping
$g_\infty:X_\infty \to \Sigma_P$ induced by $\{g_n:X_n \to S^1 \}$;
3) for each $n\in \Bbb N$ there is a point $x_n\in X_n$ such that
$g_n^{n+1}(x_{n+1})=x_n$ and $g_n(x_n)=1$. Using the classical
covering space theory ( see, e.g., \cite{13, Chapter V, Corollary
6.4}), one constructs a sequence $\{\, \phi_n: \Bbb S^1 \to X_n \,\}$  of
homeomorphisms such that $\phi_n(1) = x_n$ and $g_n \circ \phi_n = h^k_n $,
where $h^k_n$ is the $k$-fold covering mapping from (2).
In view of uniqueness of liftings \cite{13, Chapter V, Lemma 3.2},
it is easy to see that $\{\, \phi_n: n\in \Bbb N \, \}$ is a mapping between
the inverse sequences (1) and $\{X_n,g_n^{n+1}\}$.
Since each mapping $\phi_n:\Bbb S^1 \to X_n$ is a homeomorphism,
the limit mapping $\phi_\infty:\Sigma_P \to X_\infty$ induced by
$\{\,\phi_n: n \in \Bbb N \,\}$ is a homeomorphism too
\cite{6, Proposition 2.5.10}. It is clear that
$g_\infty \circ \phi_\infty = h_P^k$. In other words, the limit mappings
$g_\infty$ and $h_P^k$ are isomorphic. Hence the covering mapping
$f:X\to \Sigma_P$, being isomorphic to $g_\infty:X_\infty \to \Sigma_P$,
is isomorphic to the limit mapping $h_P^k:\Sigma_P\to \Sigma_P$ as well.

\endremark

\head
3. Periodic points
\endhead

Let $X$ be a space and let $f:X\to X$ be a mapping. For $n\in \Bbb N$,
the composite $f\circ f \circ \ldots \circ f  $ ($n$ times) is called
{\it the} $n$-{\it th iteration} of $f$ and  is denoted by $f^n$ ( $f_1=f$).
A point $x\in X$ is said to be {\it periodic} if there exists $n\in \Bbb N$
such that $f^n(x)=x$. In this case, $n$ is called {\it the period} of $x$
under $f$.

In this section we shall prove three propositions which  correspond to
the following cases:
1) $S(P)=\emptyset$; \ \ 2) $S(P)$ is an infinite set;
\ \ 3) $S(P)$ is a nonempty finite set.

We first remark that for any sequence of prime numbers $P$ the limit
mapping $h^1_P$ is the identity on the $P$-adic solenoid
$\Sigma_P$. Therefore, the set of periodic points of
$h^1_P$ coincides with the whole space $\Sigma_P$.

\proclaim{Proposition 5}
Let $P=(p_1,p_2,\ldots)$ be a sequence of prime numbers such that each prime
number from $\Bbb N$ occurs infinitely often in $P$. Then, for each
$k \geq 2$, the identity $e$ of the $P$-adic solenoid is the only
periodic point of the limit mapping $h^k_P.$
\endproclaim

\demo{Proof} Fix $k \geq 2$ and suppose that
$z=(z_1,z_2, \ldots)$ is a point of $\Sigma_P$
such that $(h^k_P)^m(z)=h^{k^m}_P(z)=z$ for some $m\in \Bbb N$.
Then, for each $n\in \Bbb N$, we have $z_n^{k^m-1}=1$ with $k^m \geq 2$.
If $k^m = 2,$ then $z=e$. Let $k^m \geq 3$.
For given $n\in \Bbb N$, we choose $l > n$  such that the product
$p_np_{n+1}\ldots p_{l-1}$ of the terms of $P$ is a multiple of $k^m-1$.
Then
$$z_n=f_n^l(z_l)=z_l^{p_np_{n+1}\ldots p_{l-1}}=1. \tag{$4$}$$
Since this is valid for each $n\in \Bbb N$, we have $z=e$, as required.
\qed
\enddemo

The next proposition is based on \cite{3, Proposition 40}
(\, cf. \cite{17, Proposition 9} \,).
Note that, in the case of the dyadic solenoid, the set $ S((2,2,\ldots))$
is infinite. Before coming to Proposition 6 we recall Euler's Theorem.
It states that for relatively prime positive integers $a$ and $m$ the
number $a^{\varphi(m)}-1$ is divisible by $m$. Here, $\varphi(m)$ is
a value of the Euler function at $m$. That is, $\varphi(1)=1$ and, for
$m>1$, $\varphi(m)=
q_1^{l_1-1}q_2^{l_2-1}\ldots q_n^{l_n-1}(q_1-1)(q_2-1)\ldots (q_n-1)$,
where $m=q_1^{l_1}q_2^{l_2}\ldots q_n^{l_n}$ is the canonical
factorization of $m$, i.e., $q_1,q_2,\ldots, q_n \geq 2$ are distinct prime
numbers and $l_1,l_2,\ldots,l_n\in \Bbb N$.

\proclaim{Proposition 6}
Let $P=(p_1,p_2, \ldots)$ be a sequence of prime numbers such that
the set $S(P)$ is infinite. Then, for each $k\in \Bbb N$, the set of
all periodic points of the limit mapping $h^k_P$ is dense in $\Sigma_P$.
\endproclaim

\demo{Proof} Given $k>1$, we choose $q\in S(P)$ such that $q > k$.
Let $\pi_n:\Sigma_P\to \Bbb S^1:(z_1,z_2,\ldots) \mapsto z_n$ be the
$n$-th projection mapping of the inverse sequence (1), $n\in \Bbb N$.

Take a basic open subset $U$ of $\Sigma_P.$ That is, $U=\pi_l^{-1}(V_l)$
for some $l\in \Bbb N$ and some open subset $V_l$ of $\Bbb S^1$.
Recall that for any $n \in \Bbb N$ such that $n\geq l$, the projection
mappings $\pi_n$ and $\pi_l$ satisfy the equality $\pi_l=f_l^n \circ \pi_n$.
Therefore, for each $n \geq l$, we have
$U=\pi_l^{-1}(V_l)=\pi_n^{-1}(V_n)$,
where $V_n=(f_l^n)^{-1}(V_l)$ is an open subset of $\Bbb S^1$.
Since $q\in S(P)$ we can choose $n\geq l$ such that $p_j\neq q$
for every $j\geq n.$ Fix such $n$.

Take $m\in \Bbb N$ such that the following condition is fulfilled:

$$\quad  \text{there exists a point} \quad z_n\in V_n \quad \text{such  that}
\quad z_n^{q^m}=1. \tag{$5$}$$

For each $j\geq n$ the numbers $p_j$ and $q^m$ are relatively prime.
Consider the automorphism
$\psi_{p_j}:\root{q^m} \of{1}\to \root{q^m} \of{1}:z\mapsto
z^{p_j}$ and its inverse $\phi_{p_j}$.
Thus we have
$$(\phi_{p_j}(z))^{p_j}=z \quad \text{for all} \quad z\in \root{q^m} \of{1}.
\tag{$6$}$$

Let $\zeta=(z_n^{p_1p_2\ldots p_{n-1}},\ldots,z_n^{p_{n-1}},z_n,
\phi_{p_n}(z_n), \phi_{p_{n+1}}\circ \phi_{p_n}(z_n), \ldots).$
According to (6), the point $\zeta$ belongs to  $\Sigma_P$. Moreover,
by (5),
$$\zeta \in U. \tag{$7$}$$

We claim that $\zeta$ is a periodic point of $h^k_P$. Indeed, since each term
of the sequence $\zeta$ is a $q^m$-th root of $1$, we get
$$\zeta^{q^m}=e.  \tag{$8$}$$

Since $k$ and $q^m$ are relatively prime, by Euler's Theorem,
$k^{q^{m-1}(q-1)}-1$ is divisible by $q^m$. Hence, by (8), \
$$\zeta^{k^{q^{m-1}(q-1)}-1}=e \qquad \text{and} \qquad
(h^k_P)^{q^{m-1}(q-1)}(\zeta)=\zeta^{k^{q^{m-1}(q-1)}}=\zeta.$$
In other words, $\zeta$ is a point of period $q^{m-1}(q-1)$ under $h^k_P$.

In view of (7) and since $U$ was chosen as an arbitrary basic open
subset of $\Sigma_P$, the proof is complete. \qed
\enddemo

\proclaim{Proposition 7}
Let $P=(p_1,p_2, \ldots)$ be a sequence of prime numbers such that
the set $S(P)$ is nonempty and finite. Let $k\in \Bbb N$ be given.
If $k$ is a multiple of the product of all prime numbers from $S(P)$,
then the identity $e$ of the $P$-adic solenoid
is the only periodic point of the limit mapping $h^k_P.$
If there exists a prime number from $S(P)$ which is not a divisor of $k$,
then the set of all periodic points of the limit mapping $h^k_P$ is dense
in $\Sigma_P$.
\endproclaim

\demo{Proof} Suppose that $S(P)=\{q_1,\ldots,q_t\}$, where $ t\in
\Bbb N$ .

First, let \ $k$ be a multiple of $q_1 \cdot \ldots \cdot q_t$,  \
and let $z=(z_1,z_2,\ldots)\in \Sigma_P$ be a periodic point of period
$m\in \Bbb N$ under $h^k_P$. Thus, for all $n\in \Bbb N$, we have
$z_n^{k^m-1}=1$, where $k^m\geq 2$. If $k^m=2$, then $z=e$.
If $k^m\geq 3$, then we choose $n\in \Bbb N$ such that
$$p_j \notin S(P)\quad \text {for all} \quad j\geq n. \tag {$9$}$$
Obviously, for each $q_j \in S(P), j\in \{1,\ldots,t\}$, the number
$k^m-1$ is not divisible by $q_j$. Hence, each prime divisor of $k^m-1$
occurs infinitely often in the sequence $P$. Choose $l > n$ such that
the number $p_np_{n+1}\ldots p_{l-1}$  is a multiple of $k^m-1$.
Then we have $z_n=1$ \ (see (4)). This implies $z_l=f_l^n(z_n)=1$
for every $l\leq n$. Since this is true for each $n$ satisfying (9), we
get  $z=e$.

Second, we assume that $k > 1$ is not divisible by $q\in S(P)$. Then every power
of $q$ and $k$ are relatively prime. Repeating the proof of
Proposition 6 for these $k$ and $q$ , one can see that the set of periodic
points of $h^k_P$ is dense in $\Sigma_P$.
\qed
\enddemo


\Refs\nofrills{References}

\ref  \no 1 \nofrills
\by J. M. Aarts and R. J. Fokkink
\paper The classification of solenoids
\jour Proc. Amer. Math. Soc. \nofrills
\vol 111 (1991)
\pages  1161 -- 1163
\endref

\ref \no 2 \nofrills
\by R. H. Bing
\paper A simple closed curve is the only homogeneous bounded plane
continuum that contains an arc
\jour Canad. J. Math. \nofrills
\vol 12 (1960)
\pages 209 -- 230
\endref

\ref \no 3
\by J. J. Charatonik and P. P. Covarrubias
\paper On covering mappings on solenoids
\jour Proc. Amer. Math. Soc.
\vol 130 (2002)
\pages 2145--2154
\endref

\ref \no 4 \nofrills
\by D. van Dantzig
\paper Ueber topologisch homogene Kontinua
\jour Fund. Math.
\vol 15 (1930)
\pages  102 -- 125
\endref

\ref \no 5
\by S. Eilenberg and N. Steenrod
\book Foundations of algebraic topology
\publ Princeton Univ. Press, Princeton, N. J.
\yr 1952
\endref

\ref \no 6
\by R. Engelking
\book General topology
\publ Monografie Matematyczne, Vol. 60, Polish Scientific Publishers, Warszawa
\yr 1977
\endref

\ref \no 7
\by R. H. Fox
\paper On shape
\jour Fund. Math.
\vol 74 (1972)
\pages  47 -- 71
\endref

\ref \no 8
\by R. H. Fox
\book Shape theory and covering spaces
\publ Lecture Notes in Math., Vol. 375, 71--90, Topology Conference Virginia
Polytechnic Institute, 1973 (R.F. Dickman Jr., P.Fletcher, eds.), Springer,
Berlin --- Heidelberg --- New York
\yr 1974
\endref

\ref \no 9
\by S. A. Grigorian and R. N. Gumerov
\paper On a covering group theorem and its applications
\jour Loba\-chevskii J. Math.
\vol X (2002)
\pages 9--16
\endref

\ref \no 10
\by S. A. Grigorian and R. N. Gumerov
\paper The structure of finite-sheeted coverings
of compact connected groups \nofrills
\jour (submitted)
\vol
\pages
\endref

\ref \no 11
\by E. Hewitt and K. A. Ross
\book Abstract harmonic analysis, Vol. I
\publ Springer-Verlag, Berlin
\yr 1963
\endref

\ref \no 12
\by S. Marde\v{s}i\'{c}, V. Matijevi\'{c}
\paper Classifying overlay structures of topological spaces
\jour Topology Appl. \nofrills
\vol 113 (2001)
\pages 167 -- 209
\endref


\ref \no 13
\by W. S. Massey
\book Algebraic topology: An Introduction
\publ Harcourt Brace Jovanovich, New York
\yr 1967
\endref

\ref \no 14
\by M. C. McCord
\paper Inverse limit sequences with covering maps
\jour Trans. Amer. Math. Soc. \nofrills
\vol 114 (1965)
\pages  197 -- 209
\endref

\ref \no 15
\by T. T. Moore
\paper On Fox's theory of overlays
\jour Fund. Math.
\vol 99 (1978)
\pages  205 -- 211
\endref


\ref \no 16
\by G. T. Whyburn
\book Analytic topology
\publ Amer. Math. Soc. Colloq. Publ. 28, Providence
\yr 1942
\endref


\ref \no 17
\by Zhou Youcheng
\paper Covering mapping on solenoids and their dynamical properties
\jour Chinese Sci. Bull.
\vol 45 (2000)
\pages  1066--1070
\endref



\endRefs
\enddocument